\documentclass{pnastwo-ppt}
\usepackage{amssymb,amsfonts,amsmath,graphicx}

\newcommand{\R}{{\mathbb R}}
\newcommand{\beq}[1]{\begin{equation}\label{#1}}
\newcommand{\eeq}{\end{equation}}
\renewcommand{\(}{\left(}
\renewcommand{\)}{\right)}

\url{}
\copyrightyear{2009}
\issuedate{Issue Date}
\volume{Volume}
\issuenumber{Issue Number}

\begin{document}

\title{Sharp rates of decay of solutions to the nonlinear fast diffusion equation via functional inequalities}

\author{M. Bonforte\affil{1}{Depto. de Matem\'{a}ticas, Universidad Aut\'{o}noma de Madrid (UAM), Campus de Cantoblanco, 28049 Madrid, Spain},
J. Dolbeault\affil{2}{Ceremade (UMR CNRS nr.~7534), Universit\'e Paris-Dau\-phine, Place de Lattre de Tassigny, 75775 Paris 16, France},
G. Grillo\affil{3}{Dip. di Matematica, Politecnico di Torino, Corso Duca degli Abruzzi 24, 10129 Torino, Italy.},
J. L. V\'azquez\affil{1}{}\hspace{-2mm}\affil{4}{ICMAT~at~UAM}
}
\maketitle

%%%%%%%%%%%%%%%%%%%%%%%%%%%%%%%%%%%%%%%%%%%%%%%%%%%%%%%%%%%%%%%%
\begin{article}

\begin{abstract} The goal of this note is to state the optimal decay rate for solutions of the nonlinear fast diffusion equation and, in self-similar variables, the optimal convergence rates to Barenblatt self-similar profiles and their generalizations. It relies on the identification of the optimal constants in some related Hardy-Poincar\'e inequalities and concludes a long series of papers devoted to generalized entropies, functional inequalities and rates for nonlinear diffusion equations.\end{abstract}


%\abbreviations{FDE, fast diffusion equation; PME, porous media equation; HPI, Hardy-Poincar\'e inequalities}

%%%%%%%%%%%%%%%%%%%%%%%%%%%%%%%%%%%%%%%%%%%%%%%%%%%%%%%%%%%%%%%%
%%%%%%%%%%%%%%%%%%%%%%%%%%%%%%%%%%%%%%%%%%%%%%%%%%%%%%%%%%%%%%%%
\section{Introduction}
The evolution equation
\beq{Eqn}
\frac{\partial u}{\partial \tau}=\nabla\cdot(u^{m-1}\,\nabla u)=\frac 1m\,\Delta u^m
\eeq
with $m\neq 1$ is a simple example of a nonlinear diffusion equation which generalizes the heat equation, and appears in a wide number of applications. Solutions differ from the linear case in many respects, notably concerning existence, regularity and large-time behaviour. We consider positive solutions $u(\tau,y)$ of this equation posed for $\tau\ge 0$ and $y\in \R^d$, $d\ge 1$. The parameter $m$ can be any real number. The equation makes sense even in the limit case $m=0$, where $u^m/m$ has to be replaced by $\log u$, and is formally parabolic for all $m\in\R$. Notice that \eqref{Eqn} is degenerate at the level $u=0$ when $m>1$ and singular when $m<1$. We consider the initial-value problem with nonnegative datum $u(\tau=0,\cdot)=u_0\in L_{{\rm loc}}^1(dx)$, where $dx$ denotes Lebesgue's measure on $\R^d$. Further assumptions on~$u_0$ are needed and will be specified later.

The description of the asymptotic behaviour of the solutions of \eqref{Eqn} as $\tau\to\infty$ is a classical and very active subject. If $m=1$, the convergence of solutions of the heat equation with $u_0\in L_+^1(dx)$ to the Gaussian kernel (up to a mass factor) is a cornerstone of the theory. In the case of equation \eqref{Eqn} with $m>1$, known in the literature as the \emph{porous medium equation}, the study of asymptotic behaviour goes back to \cite{MR586735}. The result extends to the exponents $m\in(m_c,1)$ with $m_c:=(d-2)/d$; see~\cite{MR1977429}. In these results the Gaussian kernel is replaced by some special self-similar solutions $U_{D, T}$ known as the \emph{Barenblatt solutions} (see~\cite{Ba52}) given by
\begin{equation}\label{baren.form1}
U_{D, T}(\tau,y):=\frac 1{R(\tau)^d} \left(D+\tfrac{1-m}{2\,d\,|m-m_c|}
\big|\tfrac{y}{R(\tau)}\big|^2\right)^{-\frac{1}{1-m}}_+
\end{equation}
whenever $m>m_c$ and $m\ne 1,$ with
\[
R(\tau):=(T+\tau)^{\frac{1}{d\,(m-m_c)}},
\]
where $T\ge 0$ and $D>0$ are free parameters. To some extent, these solutions play the role of the fundamental solution of the linear diffusion equations, since $\lim_{\tau\to 0}U_{D,0}(\tau,y)=M\,\delta$, where $\delta$ is the Dirac delta distribution, and $M$ depends on~$D$. Notice that the Barenblatt solutions converge as $m\to 1$ to the fundamental solution of the heat equation, up to the mass factor $M$. The results of \cite{MR586735,MR1977429} say that $U_{D,T}$ also describes the large time asymptotics of the solutions of equation \eqref{Eqn} as $\tau\to\infty$ provided $M=\int_{\R^d}u_0\,dy$ is finite, a condition that uniquely determines $D=D(M)$. Notice that in the range $m\ge m_c$, solutions of~\eqref{Eqn} with $u_0\in L_+^1(dx)$  exist globally  in time and mass is conserved: $\int_{\R^d}u(\tau,y)\,dy=M$ for any $\tau\ge 0$.

On the other hand, when $m<m_c$, a natural extension for the Barenblatt functions can be achieved by considering the same expression~\eqref{baren.form1}, but a different form for~$R$, that is
\begin{equation*}
R(\tau):=(T-\tau)^{-\frac{1}{d\,(m_c-m)}}\,.
\end{equation*}
The parameter $T$ now denotes the \emph{extinction time}, a new and important feature. The limit case $m=m_c$ is covered by $R(\tau)=e^\tau$, $U_{D,T}(\tau,y)=e^{-d\,\tau}\left(D+e^{-2\tau}\,|y|^2/d\right)^{-d/2}$. See \cite{BBDGV, MR2282669} for more detailed considerations.

In this note, we shall focus our attention on the case $m<1$ which has been much less studied. In this regime, \eqref{Eqn} is known as the \emph{fast diffusion} equation. We do not even need to assume $m>0$. We shall summarize
and extend a series of recent results on the basin of attraction of the family of \emph{generalized Barenblatt solutions} and establish the optimal rates of convergence of the solutions of \eqref{Eqn} towards a unique attracting limit state in that family. To state such a result, it is more convenient to rescale the flow and rewrite \eqref{Eqn} in self-similar variables by introducing for $m\neq m_c$ the time-dependent change of variables
\begin{equation}\label{eq:chgvariable}
t:=\tfrac{1-m}{2}\log\left(\frac{R(\tau)}{R(0)}\right)\quad\mbox{and}\quad x:=\sqrt{\tfrac{1-m}{2\,d\,|m-m_c|}}\,\frac y{R(\tau)}\,,
\end{equation}
with $R$ as above. If $m=m_c$, we take $t=\tau/d$ and $x=e^{-\tau}\,y/\sqrt d$. In these new variables, the generalized Barenblatt functions $U_{D, T}(\tau,y)$ are transformed into \emph{generalized Barenblatt profiles} $V_D(x)$, which are stationary:
\begin{equation}\label{newBaren}
V_D(x):=\(D+|x|^2\)^\frac 1{m-1}\quad x\in \R^d\,.
\end{equation}
If $u$ is a solution to~\eqref{Eqn}, the function
\begin{equation*}
v(t,x):= R(\tau)^{d}\,u(\tau,y)
\end{equation*}
solves the equation
\begin{equation}\label{FPeqn}
\frac{\partial v}{\partial t}=\nabla\cdot\left[v\,\nabla\left(\frac{v^{m-1}-V_D^{m-1}}{m-1}\right)\right]\quad t> 0\,,\quad x\in \R^d\,,
\end{equation}
with initial condition $v(t=0,x)=v_0(x):=R(0)^{-d}\,u_0(y)$ where $x$ and~$y$ are related according to \eqref{eq:chgvariable} with $\tau=0$. This nonlinear Fokker-Planck equation can also be written as
\[
\frac{\partial v}{\partial t}=\frac 1m\,\Delta v^m+\frac 2{1-m}\,\nabla\cdot(x\,v)\quad t> 0\,,\quad x\in \R^d\,.
\]

%%%%%%%%%%%%%%%%%%%%%%%%%%%%%%%%%%%%%%%%%%%%%%%%%%%%%%%%%%%%%%%%
\section{Main results}

Our main result is concerned with the \emph{sharp rate} at which a solution $v$ of the rescaled equation \eqref{FPeqn} converges to the \emph{generalized Barenblatt profile} $V_D$ given by formula \eqref{newBaren}
in the whole range $m< 1$. Convergence is measured in terms of the relative entropy given by the formula
\[
\mathcal E[v]:=\frac 1{m-1}\int_{\R^d}\left[\frac{v^m-V_D^m}{m}-\,V_D^{m-1}(v-V_D)\right]\,dx
\]
for all $m\ne 0$ (modified as mentioned for $m=0$). In order to get such convergence we need the following assumptions on the initial datum $v_0$ associated to \eqref{FPeqn}:\smallskip

\noindent {\bf (H1)} \emph{$V_{D_0}\le v_0 \le V_{D_1}$ for some $D_0>D_1>0$,}

\noindent {\bf (H2)} \emph{if $d\ge 3$ and $m\le m_*$, $(v_0-V_D)$ is integrable for a suitable $D\in[D_1,D_0]$.}

\smallskip\noindent The case $m=m_*:=(d-4)/(d-2)$ will be discussed later. Besides, if $m>m_*$, we define $D$ as the unique value in $[D_1,D_0]$ such that $\int_{{\mathbb R}^d}(v_0-V_D)\,dx=0$.
%---------------------------------------------------------------
\smallskip\par\begin{theorem}\label{Thm:Main1} Under the above assumptions, if $m<1$ and $m\ne m_*$, the entropy decays according to
\begin{equation}
\mathcal E[v(t,\cdot)]\le C\, e^{-2\,\Lambda\, t}\quad\forall\;t\ge 0\,.
\end{equation}
The sharp decay rate $\Lambda$ is equal to the best constant $\Lambda_{\alpha,d}>0$ in the Hardy--Poincar\'e inequality of Theorem~{\rm \ref{Thm:Main2}} with $\alpha:=1/(m-1)<0$. Moreover, the constant $C>0$ depends only on $m,d,D_0,D_1,D$ and $\mathcal E[v_0]$.
\end{theorem}\smallskip\par\indent
%---------------------------------------------------------------
The precise meaning of what \emph{sharp rate} means will be discussed at the end of this paper. As in \cite{BBDGV}, we can deduce from Theorem \ref{Thm:Main1} rates of convergence in more standard norms, namely, in $L^q(dx)$ for $q\ge \max\{1,d\,(1-m)/\,[2\,(2-m)+d\,(1-m)]\}$, or in $C^k$ by interpolation. Moreover, by undoing the time-dependent change of variables~\eqref{eq:chgvariable}, we can also deduce results on the \emph{intermediate asymptotics} for the solution of equation \eqref{Eqn}; to be precise, we can get rates of decay of $u(\tau,y)-R(\tau)^{-d}\,U_{D,T}(\tau,y)$ as $\tau\to+\infty$ if $m\in [m_c,1)$, or as $\tau\to T$ if $m\in(-\infty,m_c)$.

It is worth spending some words on the basin of attraction of the Barenblatt solutions $U_{D,T}$ given by~\eqref{baren.form1}. Such profiles have two parameters: $D$ corresponds to the \emph{mass} while $T$ has the meaning of the \emph{extinction time} of the solution for $m<m_c$ and of a time-delay parameter otherwise. Fix $T$ and $D$, and consider first the case $m_*< m<1$. The basin of attraction of $U_{D,T}$ contains all solutions corresponding to data which are trapped between two Barenblatt profiles $U_{D_0,T}(0,\cdot), U_{D_1,T}(0,\cdot)$ for the same value of $T$ and such that $\int_{{\mathbb R}^d}[u_0-U_{D,T}(0,\cdot)]\,dy=0$ for some $D\in [D_1,D_0]$. If $m<m_*$ the basin of attraction of a Barenblatt solution contains all solutions corresponding to data which, besides being trapped between $U_{D_0,T}$ and $U_{D_1,T}$, are integrable perturbations of $U_{D,T}(0,\cdot)$.

Now, let us give an idea of the proof of Theorem~\ref{Thm:Main1}. First assume that $D=1$ (this entails no loss of generality). On~$\R^d$, we shall therefore consider the measure $d\mu_\alpha:=h_\alpha\,dx$, where the weight $h_\alpha$ is the Barenblatt profile, defined by $h_\alpha(x):= (1+|x|^2)^{\alpha}$, with $\alpha=1/(m-1)<0$, and study on the weighted space $L^2(d\mu_\alpha)$ the operator
\[
\mathcal L_{\alpha,d}:=-h_{1-\alpha}\,\mathrm{div}\left[\,h_\alpha\,\nabla\cdot\,\right]
\]
which is such that $\int_{\R^d}f\,(\mathcal L_{\alpha,d}\,f)\,d\mu_{\alpha-1}=\int_{\R^d}|\nabla f|^2\,d\mu_\alpha$. This operator appears in the linearization of~\eqref{FPeqn} if, at a formal level, we expand
\[
v(t,x)=h_\alpha(x)\left[1+\varepsilon\,f\!\(t,x\)h_\alpha^{1-m}(x)\right]
\]
in terms of $\varepsilon$, small, and only keep the first order terms:
\[
\frac{\partial f}{\partial t}+\mathcal L_{\alpha,d}\,f=0\,.
\]
\noindent The convergence result of Theorem~\ref{Thm:Main1} follows from the energy analysis of this equation based on the Hardy-Poincar\'e inequalities that are described below. Let us fix some notations. For $d\ge 3$, let us define $\alpha_*:=-(d-2)/2$ corresponding to $m=m_*$; two other exponents will appear in the analysis, namely, $m_1:=(d-1)/d$ with corresponding $\alpha_1=-d$, and $m_2:=d/(d+2)$ with corresponding $\alpha_2=-(d+2)/2$. We have $m_*<m_c<m_2<m_1<1$. Similar definitions for $d=2$ give $m_*=-\infty$ so that $\alpha_*=0$, as well as $m_c=0$, and $m_1=m_2=1/2$. For the convenience of the reader, a table summarizing the key values of the parameter $m$ and the corresponding values of $\alpha$ is given in the Appendix.

%---------------------------------------------------------------
\smallskip\par\begin{theorem}[Sharp Hardy-Poincar\'e inequalities]\label{Thm:Main2} Let $d\ge 3$. For any $\alpha\in(-\infty,0)\setminus\{\alpha_*\}$, there is a positive constant ${\Lambda_{\alpha,d}}$
such that
\begin{equation}\label{gap}
{\Lambda_{\alpha,d}}\int_{\R^d}|f|^2\,d\mu_{\alpha-1}\leq \int_{\R^d}|\nabla f|^2\,d\mu_\alpha\quad\forall\;f\in H^1(d\mu_\alpha)
\end{equation}
under the additional condition $\int_{\R^d}f\,d\mu_{\alpha-1}=0$ if $\alpha<\alpha_*$. Moreover, the sharp constant $\Lambda_{\alpha,d}$ is given by
\[
\Lambda_{\alpha,d}=\left\{\begin{array}{ll}
\frac 14\,(d-2+2\,\alpha)^2&\mbox{if}\;\alpha\in\left[-\frac{d+2}{2},\alpha_*\right)\cup(\alpha_*,0)\,,\vspace{6pt}\cr
-\,4\,\alpha-2\,d&\mbox{if}\;\alpha\in\left[-d,-\frac{d+2}{2}\right)\,,\vspace{6pt}\cr
-\,2\,\alpha&\mbox{if}\;\alpha\in(-\infty,-d)\,.\cr
\end{array}\right.
\]
For $d=2$, inequality \eqref{gap} holds for all $\alpha<0$, with the corresponding values of the best constant $\Lambda_{\alpha,2}=\alpha^2$ for $\alpha\in [-2,0)$
and $\Lambda_{\alpha,2}=-2\alpha$ for $\alpha\in (-\infty,-2)$. For $d=1$, \eqref{gap} holds, but the values of $\Lambda_{\alpha,1}$ are given by $\Lambda_{\alpha,1}=-2\,\alpha$ if $\alpha<-1/2$ and $\Lambda_{\alpha,1}=(\alpha-1/2)^2$ if $\alpha\in[-1/2,0)$.
\end{theorem}\par\indent
%---------------------------------------------------------------

The Hardy-Poincar\'e inequalities \eqref{gap} share many properties with Hardy's inequalities, because of homogeneity reasons. A simple scaling argument indeed shows that
\[
{\Lambda_{\alpha,d}} \int_{\R^d}|f|^2\,(D+|x|^2)^{\alpha-1}\,dx\leq \int_{\R^d}|\nabla f|^2\,(D+|x|^2)^{\alpha}\,dx
\]
holds for any $f\in H^1((D+|x|^2)^{\alpha}dx)$ and any $D\ge 0$, under the additional conditions $\int_{\R^d}f\,(D+|x|^2)^{\alpha-1}\,dx=0$ and $D>0$ if $\alpha<\alpha_*$. In other words, the optimal constant, $\Lambda_{\alpha,d}$, does not depend on $D>0$ and the assumption $D=1$ can be dropped without consequences. In the limit $D\to 0$, they yield  weighted Hardy type inequalities, cf. \cite{CKN,hardy1934inequalities}.

Theorem \ref{Thm:Main2} has been proved in \cite{BBDGV} for $m<m_*$. The main improvement of this note compared to \cite{BBDGV-CRAS,BBDGV} is that we are able to give the value of the sharp constants also in the range $(m_*,1)$. These constants are deduced from the spectrum of the operator $\mathcal L_{\alpha,d}$, that we shall study below.

It is relatively easy to obtain the classical decay rates of the linear case in the limit $m\to 1$ by a careful rescaling such that weights become proportional to powers of the modified expression $(1+(1-m)\,|x|^2)^{-1/(1-m)}$. In the limit case, we obtain the Poincar\'e inequality for the Gaussian weight. As for the evolution equation, the time also has to be rescaled by a factor $(1-m)$. We leave the details to the reader. See \cite{BBDE} for further considerations on associated functional inequalities.

%%%%%%%%%%%%%%%%%%%%%%%%%%%%%%%%%%%%%%%%%%%%%%%%%%%%%%%%%%%%%%%%

\section{A brief historical overview}

The search for \emph{sharp decay rates} in fast diffusion equations has been extremely active over the last three decades. Once plain convergence of the suitably rescaled flow towards an asymptotic profile is established (cf. \cite{MR586735,MR1977429} for $m>m_c$ and \cite{BBDGV,Daskalopoulos-Sesum2006} for $m\le m_c$), getting the rates is the next step in the asymptotic analysis. An important progress was achieved by M. Del Pino and J. Dolbeault in \cite{MR1940370} by identifying sharp rates of decay for the relative entropy, that had been introduced earlier by J. Ralston and W.I. Newman in~\cite{MR760591,MR760592}. The analysis in \cite{MR1940370} uses the optimal constants in Gagliardo-Nirenberg inequalities, and these constants are computed. J.A. Carrillo and G. Toscani in \cite{MR1777035} gave a proof of decay based on the entropy/entropy-production method of D. Bakry and M. Emery, and established an analogue of the Csisz\'ar-Kullback inequality which allows to control the convergence in $L^1(dx)$, in case $m>1$. F. Otto then made the link with gradient flows with respect to the Wasserstein distance, see \cite{MR1842429}, and D. Cordero-Erausquin, B. Nazaret and C.~Villani gave a proof of Gagliardo-Nirenberg inequalities using mass transportation techniques in \cite{MR2032031}.

The condition $m\ge(d-1)/d=:m_1$ was definitely a strong limitation to these first approaches, except maybe for the entropy/entropy-production method. Gagliardo-Nirenberg inequalities degenerate into a critical Sobolev inequality for $m=m_1$, while the displacement convexity condition requires $m\ge m_1$. It was a puzzling question to understand what was going on in the range $m_c<m<m_1$, and this has been the subject of many contributions. Since one is interested in understanding the convergence towards Barenblatt profiles, a key issue is the integrability of these profiles and their moments, in terms of $m$. To work with Wasserstein's distance, it is crucial to have second moments bounded, which amounts to request $m>d/(d+2)=:m_2$ for the Barenblatt profiles. The contribution of J. Denzler and R. McCann in \cite{MR1982656,MR2126633} enters in this context. Another, weaker, limitation appears when one only requires the integrability of the Barenblatt profiles, namely $m>m_c$. Notice that the range $[m_c,1)$ is also the range for which $L^1(dx)$ initial data give rise to solutions which preserve the mass and globally exist, see for instance~\cite{MR586735,MR2282669}.

It was therefore natural to investigate the range $m\in(m_c,1)$ with entropy estimates. This has been done first by linearizing around the Barenblatt profiles in \cite{MR1901093,MR1974458}, and then a full proof for the nonlinear flow was done by J.A. Carrillo and J.L. V\'azquez in \cite{MR1986060}. A detailed account for these contributions and their motivations can be found in the survey paper \cite{MR2065020}. Compared to classical approaches based on comparison, as in the book~\cite{MR2282669}, a major advantage of entropy techniques is that they combine very well with $L^1(dx)$ estimates if $m>m_c$, or relative mass estimates otherwise, see~\cite{BBDGV}.

 The picture for $m\le m_c$  turns out to be entirely different and more complicated, and it was not considered until quite recently. First of all, many classes of solutions vanish in finite time, which is a striking property that forces us to change the concept of asymptotic behaviour from large-time behaviour to behaviour near the extinction time. On the other hand, $L^1(dx)$ solutions lose mass as time evolves. Moreover, the natural extensions of Barenblatt's profiles make sense but these profiles have two novel properties: they vanish in finite time and they do not have finite~mass.

There is a large variety of possible behaviours and many results have been achieved, such as the ones described in \cite{MR2282669} for data which decay strongly as $|x|\to\infty$. However, as long as one is interested in solutions converging towards Barenblatt profiles in self-similar variables, there were some recent  results on plain convergence: a paper of P. Daskalopoulos and N. Sesum, \cite{Daskalopoulos-Sesum2006}, using comparison techniques, and in two contributions involving the authors of this note, using \emph{relative entropy methods}, see \cite{BBDGV-CRAS,BBDGV}. This last approach proceeds further into the description of the convergence by identifying a suitable weighted linearization of the relative entropy. In the appropriate space, $L^2(d\mu_{\alpha-1})$, with the notations of Theorem~\ref{Thm:Main2}, it gives rise to an exponential convergence after rescaling. This justifies the heuristic computation which relates Theorems~\ref{Thm:Main1} and \ref{Thm:Main2}, and allows to identify the sharp rates of convergence. The point of this note is to explicitly state and prove such  rates in the whole range $m<1$.

%%%%%%%%%%%%%%%%%%%%%%%%%%%%%%%%%%%%%%%%%%%%%%%%%%%%%%%%%%%%%%%%
\section{Relative entropy and linearization}

The strategy developed in \cite{BBDGV} is based on the extension of the \emph{relative entropy} of J. Ralston and W.I. Newmann, which can be written in terms of $w=v/V_D$ as
\[
\mathcal F[w]:=\frac 1{1-m}\int_{\R^d}\left[w-1-\frac{1}{m}\big(w^m-1\big)\right]\,V_D^m\,dx\,.
\]
For simplicity, assume $m\neq 0$. Notice that ${\mathcal F}[w]={\mathcal E}[v]$. Let
\[
\mathcal I[w]:=\int_{\R^d}\left|\frac 1{m-1}\,\nabla\left[(w^{m-1}-1)\,V_D^{m-1}\right]\,\right|^2v\,dx
\]
be the \emph{generalized relative Fisher information}. If $v$ is a solution of \eqref{FPeqn}, then
\beq{Eqn:Entropy-EntropyProduction}
\frac d{dt}\mathcal F[w(t,\cdot)]=-\,\mathcal I [w(t,\cdot)]\quad\forall\;t> 0
\ee
and, as a consequence, $\lim_{t\to+\infty}\mathcal F[w(t,\cdot)]=0$ for all $m<1$. The method is based on Theorem~\ref{Thm:Main2} and uniform estimates that relate linear and nonlinear quantities. Following \cite{BBDGV,BGV} we can first estimate from below and above the entropy $\mathcal{F}$ in terms of its linearization, which appears in~ \eqref{gap}:
\begin{equation}\label{Entr.lin.nonlin}
h^{m-2}\int_{\R^d}\kern -5pt|f|^2\,V_D^{2-m}\;dx\le 2\,\mathcal F[w]\le h^{2-m}\int_{\R^d}\kern -5pt|f|^2\,V_D^{2-m}\;dx
\end{equation}
where $f:=(w-1)\,V_D^{m-1}$, $h_1(t):=\mathrm{inf}_{\R^d}w(t,\cdot)$, $h_2(t):=\mathrm{sup}_{\R^d}w(t,\cdot)$ and $h:=\max\{h_2,1/h_1\}$. We notice that \hbox{$h(t)\to 1$} as $t\to+\infty$. Similarly, the generalized Fisher information satisfies the bounds
\begin{equation}\label{Fish.lin.nonlin}
\int_{\R^d}|\nabla f|^2\,V_D\;dx\le [1+X(h)]\,\mathcal I[w]+Y(h)\int_{\R^d}|f|^2\,V_D^{2-m}\;dx
\end{equation}
where $h_2^{2(2-m)}/h_1\le h^{5-2m}=:1+X(h)$ and \hbox{$d\,(1-m)$} $\big[\left(h_2/h_1\right)^{2(2-m)}-1\big]\le d\,(1-m)\,\big[h^{4(2-m)}-1\big]=:Y(h)$. Notice that $X(1)=Y(1)=0$. Joining these inequalities with the Hardy-Poincar\'e inequality of Theorem~\ref{Thm:Main2} gives
\begin{equation}
\mathcal F[w]\le \frac{h^{2-m}\,[1+X(h)]}{2\,\big[\Lambda_{\alpha,d}-Y(h)\big]}\,\mathcal{I}[w]
\end{equation}
as soon as $0<h<h_*:=\min\{h>0\,:\,\Lambda_{\alpha,d}-Y(h)\ge 0\}$. On the other hand, uniform relative estimates hold, according to \cite{BGV}, formula (5.33): for some $\mathsf C=\mathsf C(d,m,D,D_0,D_1)$,
\begin{equation}\label{Unif}
0\le h-1\le\mathsf C\,\mathcal{F}^\frac{1-m}{d+2-(d+1)m}\,.
\end {equation}
Summarizing, we end up with a system of nonlinear differential inequalities, with $h$ as above and, at least for any $t>t_*$, $t_*>0$ large enough,
\beq{Ineq:Gronwall}
\frac{d}{dt}\mathcal{F}[w(t,\cdot)]
\le-2\,\frac{\Lambda_{\alpha,d}-Y(h)}{\big[1+X(h)\big]\,h^{2-m}}\,\mathcal{F}[w(t,\cdot)]\,.
\eeq
Gronwall type estimates then show that
\[
\limsup_{t\to\infty}\,{\rm e}^{2\,\Lambda_{\alpha,d}\,t}\mathcal{F}[w(t,\cdot)]<+\infty\,.
\]
This completes the proof of Theorem~\ref{Thm:Main1} for $m\neq 0$. The adaptation to the logarithmic nonlinearity is left to the reader. Results in \cite{BBDGV} are improved in two ways: a time-dependent estimate of $h$ is used in place of $h(0)$, and the precise expression of the rate is established. One can actually get a slightly more precise estimate by coupling \eqref{Unif} and~\eqref{Ineq:Gronwall}.
%---------------------------------------------------------------
\smallskip\par\noindent\begin{corollary}\label{Cor:EntropyAndUniformEstimates}Under the assumptions of Theorem~{\rm \ref{Thm:Main1}}, if $h(0)<h_*$, then $\mathcal{F}[w(t,\cdot)]\le G\big(t,h(0),\mathcal{F}[w(0,\cdot)]\big)$ for any $t\ge 0$, where $G$ is the unique solution of the nonlinear ODE
\[
\frac{dG}{dt}
 =-2\frac{\Lambda_{\alpha,d}-Y(h)}{[1+X(h)]\,h^{2-m}}\,G
 \quad\mbox{with}\quad h=1+\mathsf C\,G^\frac{1-m}{d+2-(d+1)m}
\]
and initial condition $G(0)=\mathcal F[w(0,\cdot)]$.
\end{corollary}
%---------------------------------------------------------------

%%%%%%%%%%%%%%%%%%%%%%%%%%%%%%%%%%%%%%%%%%%%%%%%%%%%%%%%%%%%%%%%
\section{Operator equivalence. The spectrum of $\mathcal L_{\alpha,d}$}

An important point of this note is the computation of the spectrum of $\mathcal{L}_{\alpha,d}$ for any $\alpha<0$. This spectrum was only partially understood in \cite{BBDGV, BBDGV-CRAS}. In particular, the existence of a spectral gap was established for all $\alpha\ne \alpha_*= (2-d)/2$, but its value was not stated for all values of $\alpha$.

J. Denzler and R.J. McCann in \cite{MR1982656,MR2126633} formally linearized the fast diffusion flow (considered as a gradient flow of the entropy with respect to the Wasserstein distance) in the framework of mass transportation, in order to guess the asymptotic behaviour of the solutions of \eqref{Eqn}. This leads to a different functional setting, with a different linearized operator, $\mathcal H_{\alpha,d}$. They performed the detailed analysis of its spectrum for all $m\in(m_c,1)$, but the justification of the nonlinear asymptotics could not be completed due to the difficulties of the functional setting, especially in the very fast diffusion range.

Our approach is based on relative entropy estimates and the Hardy-Poincar\'e inequalities of Theorem~\ref{Thm:Main2}. The asymptotics have readily been justified in \cite{BBDGV}. The operator $\mathcal L_{\alpha,d}=-h_{1-\alpha}\,\mathrm{div}\,\left[\,h_\alpha\,\nabla\cdot\,\right]$ can be initially defined on $\mathcal D(\mathbb{R}^d)$. To construct a self-adjoint extension of such an operator, one can consider the quadratic form $f\mapsto(f,\mathcal L_{\alpha,d}\,f)$, where $(\,\cdot\,,\,\cdot\,)$ denotes the scalar product on $L^2(d\mu_{\alpha-1})$. Standard results show that such a quadratic form is closable, so that its closure defines a unique self-adjoint operator, its Friedrich's extension, still denoted by the same symbol for brevity. The operator $\mathcal H_{\alpha,d}$ is different: it is obtained by taking the operator closure of $\mathcal L_{\alpha,d}$, initially defined on $\mathcal D(\mathbb{R}^d)$, in the Hilbert space \emph{$H^{1,*}(d\mu_{\alpha}):=\big\{f\in L^2(d\mu_{\alpha-1}):$ \hbox{$\int_{\R^d}|\nabla f|^2\,d\mu_\alpha<\infty$} and \hbox{$\int_{\R^d}f\,d\mu_{\alpha-1}=0$} if $\alpha<\alpha_*\big\}$}, so that
\[
\mathcal H_{\alpha,d}\,f:=h_{1-\alpha}\,\nabla\cdot\big[h_\alpha\,\nabla \big(h_{1-\alpha}\,\nabla\cdot(h_\alpha\,\nabla f) \big) \big]\,.
\]
Because of the Hardy-Poin\-car\'e inequality, $\int_{\R^d}|\nabla f|^2\,d\mu_\alpha$ defines a~norm. Denote by $\langle\,\cdot\,,\,\cdot\,\rangle$ the corresponding scalar product and notice that $\langle\,f\,,\,g\,\rangle\!=\!(\,f\,,\,\mathcal L_{\alpha,d}\,g\,)$ for any \hbox{$f,g\in\mathcal D(\R^d)$}.
%---------------------------------------------------------------
\smallskip\par\begin{proposition}\label{Prop:Equivalence} The operator $\mathcal L_{\alpha,d}$ on $L^2(d\mu_{\alpha-1})$ has the same spectrum as the operators $\mathcal H_{\alpha,d}$ on $H^{1,*}(d\mu_{\alpha})$.
\end{proposition}\smallskip\par\indent
%---------------------------------------------------------------
The proof is based on the construction of a suitable unitary operator $U:H^{1,*}(d\mu_{\alpha})\to L^2(d\mu_{\alpha-1})$, such that $U\,\mathcal H_{\alpha,d}\,U^{-1}=\mathcal L_{\alpha,d}$. We claim that $U=\sqrt{\mathcal L_{\alpha,d}}$ is the requested unitary operator. By definition, $\mathcal D(\mathbb{R}^d)$ is a form core of $\mathcal L_{\alpha,d}$, and as a consequence, the identity has to be established only for functions $f\in\mathcal D(\mathbb{R}^d)$. Since $\|Uf\|^2=\big(\,f\,,\,\mathcal L_{\alpha,d} f\,\big)=\int_{\R^d}|\nabla f|^2\,d\mu_\alpha$, we get
\[
\begin{split}
&\hspace*{-8pt}\big(\,U\,\mathcal H_{\alpha,d}\,U^{-1}f\,,\, g\,\big)
=\big\langle\,\mathcal H_{\alpha,d}\,U^{-1}f\,,\, U^{-1}g\,\big\rangle\\
&=\big(\,\mathcal L_{\alpha,d}\,U^{-1}f\,,\, \mathcal L_{\alpha,d}\,U^{-1}g\,\big)
=\big(\,U\,f\,,\,U\,g\,\big)
=\big\langle\, f\,,\,g\,\big\rangle\\
\end{split}
\]
where we have used the properties $U^*=U^{-1}$ and $U^2=\mathcal L_{\alpha,d}$. This unitary equivalence between $\mathcal L_{\alpha,d}$ and $\mathcal H_{\alpha,d}$ implies the identity of their spectra.

\medskip We may now proceed with the presentation of the actual values of the spectrum by extending the results of \cite{MR2126633}. According to \cite{MR0282313}, the spectrum of the Laplace-Beltrami operator on $S^{d-1}$ is described by
\[
-\Delta_{S^{d-1}}Y_{\ell \mu}=\ell\,(\ell+d-2)\,Y_{\ell \mu}
\]
with $\ell=0$, $1$, $2$, \ldots and $\mu=1$, $2$, \ldots $M_\ell:=\frac{(d+\ell-3)!\,(d+2\ell-2)}{\ell!\,(d-2)!}$ with the convention $M_0=1$, and $M_1=1$ if $d=1$. Using spherical coordinates and separation of variables, the discrete spectrum of $\mathcal L_{\alpha,d}$ is therefore made of the values of $\lambda$ for which
\begin{equation}\label{EigenF}
v''+\(\tfrac{d-1}r+\tfrac{2\,\alpha\,r}{1+r^2}\)\,v'+\(\tfrac\lambda{1+r^2}-\tfrac{\ell\,(\ell+d-2)}{r^2}\)\,v=0
\end{equation}
has a solution on $\R^+\ni r$, in the domain of $\mathcal L_{\alpha,d}$. The change of variables $v(r)=r^\ell\,w(-r^2)$ allows to express $w$ in terms of the hypergeometric function $_2F_1(a,b,c;z)$ with $c=\ell+d/2$, $a+b+1=\ell+\alpha+d/2$ and $a\,b=(2\,\ell\,\alpha+\lambda)/4$, as the solution for $s=-r^2$ of
\[
s\,(1-s)\,y''+[\,c-(a+b+1)\,s\,]\,y'-a\,b\,y=0\,,
\]
see \cite{weisstein2005hf}. Based on \cite{BBDGV-CRAS,MR2126633}, we can state the following result.
%---------------------------------------------------------------
\smallskip\par\begin{proposition}\label{Prop:Spectrum} The bottom of the continuous spectrum of the operator $\mathcal L_{\alpha,d}$ on $L^2(d\mu_{\alpha-1})$ is $\lambda_{\alpha,d}^{\rm cont}:=\frac 14(d+2\,\alpha-2)^2$. Moreover, $\mathcal L_{\alpha,d}$ has some discrete spectrum only for $m>m_2=d/(d+2)$. For $d\ge 2$, the discrete spectrum is made of the eigenvalues
\begin{equation}\label{eigen}
\lambda_{\ell k}=-2\,\alpha\,\(\ell+2\,k\)-4\,k\,\(k+\ell+\frac d2-1\)
\end{equation}
with $\ell$, $k=0$, $1$, \ldots provided $(\ell,k)\neq(0,0)$ and $\ell+2k-1<-(d+2\,\alpha)/2$. If $d=1$, the discrete spectrum is made of the eigenvalues $\lambda_k=k\,(1-2\,\alpha-k)$ with $k\in{\mathbb N}\cap[1,1/2-\alpha]$.\end{proposition}\smallskip\par\indent
%---------------------------------------------------------------

\begin{figure}\begin{center}\includegraphics{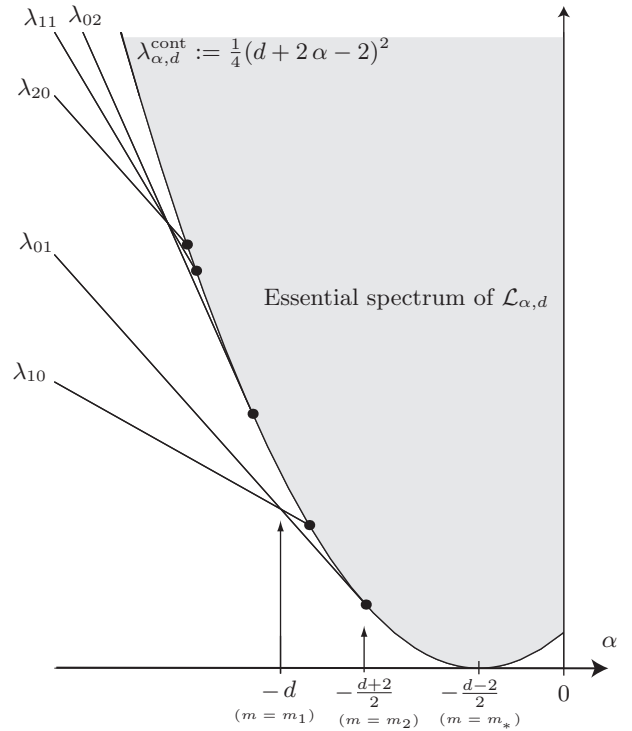}\caption{Spectrum of $\mathcal L_{\alpha,d}$ as a function of $\alpha$, for $d=5$.}\end{center}\end{figure}

Using Persson's characterization of the continuous spectrum, see \cite{MR0133586,BBDGV-CRAS}, one can indeed prove that $\lambda_{\alpha,d}^{\rm cont}$ is the optimal constant in the following inequality: for any $f\in\mathcal D(\R^d\setminus\{0\})$,
\[
\lambda_{\alpha,d}^{\rm cont}\int_{\R^d}|f|^2\,|x|^{2(\alpha-1)}\;dx\le\int_{\R^d}|\nabla f|^2\,|x|^{2\,\alpha}dx\,.
\]
The condition that the solution of \eqref{EigenF} is in the domain of $\mathcal L_{\alpha,d}$ determines the eigenvalues. A more complete discussion of this topic can be found in \cite{MR2126633}, which justifies the expression of the discrete spectrum.

\noindent Since $\alpha=1/(m-1)$, we may notice that for $d\ge 2$, $\alpha=-d$ (corresponding to $-2\,\alpha=\lambda_{10}=\lambda_{01}=-4\,\alpha-2\,d$) and $\alpha=-(d+2)/2$ (corresponding to $\lambda_{01}=\lambda_0^{\rm cont}:=\frac 14\,(d+2\,\alpha-2)^2$) respectively mean $m=m_1=(d-1)/d$ and $m=m_2=d/(d+2)$.

The above spectral results hold exactly in the same form when $d=2$, see \cite{MR2126633}. Notice in particular that $\lambda^{{\rm cont}}_{\alpha,d=2}=\alpha^2$ so that there is no equivalent of $m_*$ for $d=2$. With the notations of Theorem \ref{Thm:Main2}, $\alpha_*=0$. All results of Theorem \ref{Thm:Main1} hold true under the sole assumption (H1).

In dimension $d=1$, the spectral results are different, see \cite{MR2126633}. The discrete spectrum is nonempty whenever $\alpha\le-1/2$, that is $m\ge-1$.

%%%%%%%%%%%%%%%%%%%%%%%%%%%%%%%%%%%%%%%%%%%%%%%%%%%%%%%%%%%%%%%%
\section{The critical case}

Since the spectral gap of $\mathcal L_{\alpha,d}$ tends to zero as $m\to m_*$\,, the previous strategy fails when $m=m_*$ and one might expect a slower decay to equilibrium, sometimes referred as \emph{slow asymptotics.} The following result has been proved in~\cite{BGV}.
%---------------------------------------------------------------
\smallskip\par\begin{theorem}\label{mstar} Assume that $d\ge 3$, let $v$ be a solution of~\eqref{FPeqn} with $m= m_*$, and suppose that {\rm (H1)-(H2)} hold. If $|v_0-V_D|$ is bounded a.e. by a radial $L^1(dx)$ function, then there exists a positive constant $C^*$ such that
\begin{equation}
\mathcal E[v(t,\cdot)]\le C^*\,t^{-1/2}\quad\forall\;t\ge 0\,,
\end{equation}
where $C^*$ depends only on $m,d,D_0,D_1,D$ and $\mathcal E[v_0]$.
\end{theorem}\smallskip\par\indent
%---------------------------------------------------------------
Rates of convergence in $L^q(dx)$, $q\in(1,\infty]$ follow. Notice that in dimension $d=3$ and $4$, we have respectively $m_*=-1$ and $m_*=0$. In the last case, Theorem \ref{mstar} applies to the logarithmic diffusion.

The proof relies on identifying first the asymptotics of the linearized evolution. In this case, the bottom of the continuous spectrum of $\mathcal{L}_{\alpha_*,d}$ is zero. This difficulty is overcome by noticing that the operator $\mathcal{L}_{\alpha_*,d}$ on $\R^d$ can be identified with the Laplace-Beltrami operator for a suitable conformally flat metric on $\R^d$, having positive Ricci curvature. Then the on-diagonal heat kernel of the linearized generator behaves like $t^{-d/2}$ for small $t$ and like $t^{-1/2}$ for large $t$. The Hardy-Poincar\'e inequality is replaced by a weighted Nash inequality: there exists a positive continuous and monotone function $\mathcal N$ on $\mathbb R^+$ such that for any nonnegative smooth function $f$ with $M=\int_{\R^d}f\,d\mu_{-d/2}$ (recall that $\alpha_*-1=-d/2$),
\[
\frac 1{M^2}\,\int_{\R^d}|f|^2\,d\mu_{-d/2}\le \mathcal N\left(\frac 1{M^2}\,\int_{\R^d}|\nabla f|^2\,\;d\mu_{(2-d)/2}\right).
\]
The function $\mathcal N$ behaves as follows: $\lim_{s\to 0^+}s^{-1/3}\,\mathcal N(s)=c_1>0$ and $\lim_{s\to\infty}s^{-d/(d+2)}\mathcal N(s)=c_2>0$. Only the first limit matters for the asymptotic behaviour. Up to technicalities, inequality \eqref{Ineq:Gronwall} is replaced by $(\mathcal{F}[w(t,\cdot)])^3\le K\,\mathcal{I}[w(t,\cdot)]$ for some $K>0$, $t\ge t_0$ large enough, which allows to complete the proof.

%%%%%%%%%%%%%%%%%%%%%%%%%%%%%%%%%%%%%%%%%%%%%%%%%%%%%%%%%%%%%%%%
\section{Faster convergence}

A very natural issue is the question of improving the rates of convergence by imposing restrictions on the initial data. Results of this nature have been observed in \cite{MR1901093} in case of radially symmetric solutions, and are carefully commented in~\cite{MR2126633}. By locating the center of mass at zero, we are able to give an answer, which amounts to kill the $\lambda_{10}$ mode, whose eigenspace is generated by $x\mapsto x_i$, $i=1$, $2$\ldots $d$. This is an improvement compared to the first result in this direction, which has been obtained by R. McCann and D. Slep\v{c}ev in~\cite{MR2211152}, since we obtain an improved \emph{sharp rate} of convergence of the solution of \eqref{FPeqn}, as a consequence of the following improved Hardy-Poincar\'e inequality.
%---------------------------------------------------------------
\smallskip\par\begin{lemma}\label{Lem:Improved2} Let $\widetilde{\Lambda}_{\alpha,d}:=-4\,\alpha-2\,d$ if $\alpha<-d$ and $\widetilde{\Lambda}_{\alpha,d}:=\lambda_{\alpha,d}^{\rm cont}$ if $\alpha\in[-d,-d/2)$. If $d\ge 2$, for any $\alpha\in(-\infty,-d)$, we have
\[
\widetilde{\Lambda}_{\alpha,d} \int_{\R^d}|f|^2\,d\mu_{\alpha-1}\leq\int_{\R^d}|\nabla f|^2\,d\mu_\alpha\quad\forall\;f\in H^1(d\mu_\alpha)
\]
under the conditions $\int_{\R^d}f\,d\mu_{\alpha-1}=0$ and $\int_{\R^d}x\,f\,d\mu_{\alpha-1}=0$. The constant $\widetilde{\Lambda}_{\alpha,d}$ is sharp.
\end{lemma}\smallskip\par\indent
%---------------------------------------------------------------
This covers the range $m\in(m_1,1)$ with \hbox{$m_1=(d-1)/d$}.
%---------------------------------------------------------------
\smallskip\par\begin{theorem}\label{Thm:Main3} Assume that $m\in(m_1,1)$, $d\ge 3$. Under Assumption {\rm (H1)}, if $v$ is a solution of \eqref{FPeqn} with initial datum~$v_0$ such that $\int_{\R^d}x\,v_0\,dx=0$ and if $D$ is chosen so that $\int_{\R^d} (v_0-V_D)\,dx=0$, then there exists a positive constant~$\widetilde{C}$ depending only on $m,d,D_0,D_1,D$ and $\mathcal E[v_0]$ such that the relative entropy decays like
\[
\mathcal E[v(t,\cdot)]\le \widetilde{C}\, {\rm e}^{-\widetilde{\Lambda}_{\alpha,d}\,t}\quad\forall\;t\ge 0\;.
\]
\end{theorem}
%---------------------------------------------------------------

%%%%%%%%%%%%%%%%%%%%%%%%%%%%%%%%%%%%%%%%%%%%%%%%%%%%%%%%%%%%%%%%
\section{A variational approach of sharpness}

Recall that $(d-2)/d=m_c<m_1=(d-1)/d$. The entropy / entropy production inequality obtained in \cite{MR1940370} in the range $m\in[m_1,1)$ can be written as $\mathcal F\leq\frac 12\,\mathcal I$ and it is known to be sharp as a consequence of the optimality case in Gagliardo-Nirenberg inequalities. Moreover, equality is achieved if and only if $v=V_D$. The inequality has been extended in \cite{MR1986060} to the range $m\in(m_c,1)$ using the Bakry-Emery method, with the same constant $1/2$, and again equality is achieved if and only if $v=V_D$, but sharpness of $1/2$ is not as straightforward for $m\in(m_c,m_1)$ as it is for $m\in[m_1,1)$. The question of the optimality of the constant can be reformulated as a variational problem, namely to identify the value of the positive constant
\[
\mathcal C=\inf\frac{\mathcal I[v]} {\mathcal E [v]}
\]
where the infimum is taken over the set of all functions such that $v\in\mathcal D(\R^d)$ and $\int_{\R^d}v\,dx=M$. Rephrasing the sharpness results, we know that $\mathcal C=2$ if $m\in(m_1,1)$ and $\mathcal C\ge 2$ if $m\in(m_c,m_1)$. By taking $v_n=V_D\,(1+\frac 1n\,f\,V_D^{1-m})$ and letting $n\to\infty$, we get
\[
\lim_{n\to\infty}\frac{\mathcal I[v_n]}{\mathcal E [v_n]}=\frac{\int_{\R^d}|\nabla f|^2\,V_D\,dx}{\int_{\R^d}|f|^2\,V_D^{2-m}\,dx}\;.
\]
With the optimal choice for $f$, the above limit is less or equal than $2$. Since we already know that $\mathcal C\ge 2$, this shows that $\mathcal C= 2$ for any $m>m_c$. It is quite enlightening to observe that optimality in the quotient gives rise to indetermination since both numerator and denominator are equal to zero when $v=V_D$. This also explains why it is the first order correction which determines the value of $\mathcal C$, and, as a consequence, why the optimal constant, $\mathcal C=2$, is determined by the linearized problem.

When $m\le m_c$, the variational approach is less clear since the problem has to be constrained by a uniform estimate. Proving that any minimizing sequence $(v_n)_{n\in\mathbb N}$ is such that $v_n/V_D-1$ converges, up to a rescaling factor, to a function~$f$ associated to the Hardy-Poincar\'e inequalities would be a significant step, except that one has to deal with compactness issues, test functions associated to the continuous spectrum and a uniform constraint.

%%%%%%%%%%%%%%%%%%%%%%%%%%%%%%%%%%%%%%%%%%%%%%%%%%%%%%%%%%%%%%%%
\section{Sharp rates of convergence and conjectures}

In Theorem~\ref{Thm:Main1}, we have obtained that the rate $\exp(-\Lambda_{\alpha,d}\,t)$ is \emph{sharp.} The precise meaning of this claim is that
\[
\Lambda_{\alpha,d}=\liminf_{h\to 0_+}\inf_{w\in\mathcal S_h}\frac{\mathcal I[w]}{\mathcal F[w]}\,,
\]
where the infimum is taken on the set $\mathcal S_h$ of smooth, nonnegative bounded functions $w$ such that \hbox{$\|w-1\|_{L^\infty(dx)}\le h$} and such that $\int_{\R^d}(w-1)\,V_D\,dx$ is zero if $d=1,2$ and $m<1$, or if $d\ge 3$ and $m_*<m<1$, and it is finite if $d\ge 3$ and $m<m_*$. Since, for a solution $v(t,x)=w(t,x)\,V_D(x)$ of \eqref{FPeqn}, \eqref {Eqn:Entropy-EntropyProduction} holds, by \emph{sharp rate} we mean the best possible rate, which is uniform in $t\ge 0$. In other words, for any $\lambda>\Lambda_{\alpha,d}$, one can find some initial datum in $\mathcal S_h$ such that the estimate $\mathcal F[w(t,\cdot)]\le\mathcal F[w(0,\cdot)]\,\exp(-\lambda\,t)$ is wrong for some $t>0$. We did not prove that the rate $\exp(-\Lambda_{\alpha,d}\,t)$ is \emph{globally sharp} in the sense that for some special initial data, $\mathcal F[w(t,\cdot)]$ decays exactly at this rate, or that $\liminf_{t\to\infty}\exp(\Lambda_{\alpha,d}\,t)\,\mathcal F[w(t,\cdot)]>0$, which is possibly less restrictive.

However, if $m\in(m_1,1)$, $m_1=(d-1)/d$, then $\exp(-\Lambda_{\alpha,d}\,t)$ is also a \emph{globally sharp} rate, in the sense that the solution with initial datum $u_0(x)=V_D(x+x_0)$ for any $x_0\in\R^d\setminus\{0\}$ is such that $\mathcal F[w(t,\cdot)]$ decays exactly like $\exp(-\Lambda_{\alpha,d}\,t)$. This formally answers the dilation-persistence conjecture as formulated in \cite{MR2126633}. The question is still open when $m\le m_1$.

Another interesting issue is to understand if improved rates, that is rates of the order of $\exp(-\lambda_{\ell k}\,t)$ with $(\ell,k)\neq(0,0)$, $(0,1)$, $(1,0)$ are \emph{sharp} or \emph{globally sharp} under additional moment-like conditions on the initial data. It is also open to decide whether $\exp(-\widetilde{\Lambda}_{\alpha,d}\,t)$ is sharp or globally sharp under the extra condition $\int_{\R^d}x\,v_0\,dx=0$.

%%%%%%%%%%%%%%%%%%%%%%%%%%%%%%%%%%%%%%%%%%%%%%%%%%%%%%%%%%%%%%%%
\section{Appendix. A table of correspondence}
For the convenience of the reader, a table of definitions of the key values of $m$ when $d\ge 3$ is provided below with the correspondence for the values of $\alpha=1/(m-1)$. This note is restricted to the case $m\in(-\infty,1)$, that is $\alpha\in(-\infty,0)$.
\begin{table}[ht]
\begin{center}{\small
\begin{tabular}{r|cccccc}\hline
$m=$ &$-\infty$&$m_*$&$m_c$&$m_2$&$m_1$&$1$\cr\hline
$m=$ &$-\infty$&$\frac{d-4}{d-2}$&$\frac{d-2}d$&$\frac d{d+2}$&$\frac{d-1}d$&$1$\cr\hline
$\alpha=$ &$0$&$-\frac{d-2}2$&$-\frac d2$&$-\frac{d+2}2$&$-d$&$-\infty$\cr\hline
\end{tabular}}
\end{center}
\end{table}

%%%%%%%%%%%%%%%%%%%%%%%%%%%%%%%%%%%%%%%%%%%%%%%%%%%%%%%%%%%%%%%%
\begin{acknowledgments}
This work has been supported by the ANR-08-BLAN-0333-01 project CBDif-Fr and the exchange program of University Paris-Dauphine and Universidad Aut\'{o}noma de Madrid. MB and JLV partially supported by Project MTM2008-06326-C02-01 (Spain). MB, GG and JLV partially supported by HI2008-0178 (Italy-Spain).

\smallskip\par

\noindent
\copyright\,2009 by the authors. This paper may be reproduced, in its entirety, for non-commercial purposes.
\end{acknowledgments}
%%%%%%%%%%%%%%%%%%%%%%%%%%%%%%%%%%%%%%%%%%%%%%%%%%%%%%%%%%%%%%%%
%%%%%%%%%%%%%%%%%%%%%%%%%%%%%%%%%%%%%%%%%%%%%%%%%%%%%%%%%%%%%%%%
%\bibliographystyle{pnas}
%\bibliography{References}
%%%%%%%%%%%%%%%%%%%%%%%%%%%%%%%%%%%%%%%%%%%%%%%%%%%%%%%%%%%%%%%%
\def\cprime{$'$}

%%%%%%%%%%%%%%%%%%%%%%%%%%%%%%%%%%%%%%%%%%%%%%%%%%%%%%%%%%%%%%%%
\end{article}
\end{document}